\documentclass{amsart}
\usepackage{amsmath}
\newtheorem{theorem}{Theorem}[section]

\usepackage{cite}
\theoremstyle{definition}
\newtheorem{definition}[theorem]{Definition}

\usepackage{graphicx}
\usepackage{epstopdf}
\theoremstyle{remark}
\newtheorem{remark}[theorem]{Remark}

\numberwithin{equation}{section}

\begin{document}

\title{The Calculus of bivariate fractal interpolation surfaces}

\author{Subhash Chandra}
\address{School of Basic Sciences, Indian Institute of Technology Mandi,\\ Kamand (H.P.) - 175005, India}
\email{schandra.iitd15@gmail.com; abbas@iitmandi.ac.in}
\author{Syed Abbas}


\begin{abstract}
In this article, we investigate partial integrals and partial derivatives of bivariate fractal interpolation functions. We prove also that the mixed Riemann-Liouville fractional integral and derivative of order $\gamma = (p,q);~ p>0, q>0$, of bivariate fractal interpolation functions are again bivariate interpolation functions corresponding to some iterated function system (IFS). Furthermore, we discuss the integral transforms and fractional order integral transforms of the bivariate fractal interpolation functions.
\end{abstract}

\maketitle


\section{INTRODUCTION}

Fractal interpolation function (FIF) was introduced by Barnsley \cite{MF1,B}  through iterated function system (IFS). FIF is an interpolation function whose graph is an invariant set of an any IFS. Wang\cite{W} introduced a new class of fractal interpolation functions (FIFs) by using variable parameters. The IFS model is described by the property of self-affinity. The self-affinity specifies that the shape is invariant regardless to any degree of magnification. The calculus of FIFs for the uinvariate case is given by Barnsley\cite{MF2}. It is natural to aries question that can we present analogous results for bivariate fractal interpolation surfaces (FISs). FISs described by Massopust \cite{PM} on triangles, where the interpolation points at the boundary must be coplanar. Geronimo and Hardin \cite{GH} have generalized this work by allowing the construction over polygonal regions by taking arbitrary interpolation points. In another work by Zhao \cite{Zhao}, the author uses triangulation to construct the fractal surfaces. Dalla in \cite{Dalla} has constructed (FISs) on rectangular grids by taking the collinear interpolation points on the boundary. In \cite{Feng}, Feng considered more general construction of FISs on the rectangular grids under some restriction. The given restriction for continuity is difficult to verify. Malysz\cite{Mal} introduced fold-out technique. Metzler and Yun \cite{Metz} have generalised this method. The reader, who are interested to know about more general constructions of fractal interpolation functions, we refer to \cite{D1,D2,Xie} and references therein. Huo-Jun Ruan and Qiang Xu\cite{Ruan} presented a more general construction using fold-out technique and deduce bilinear FIS. Verma and Viswanathan\cite{V1,V2} has discussed fractal operator which is associated with bivariate FIFs and about the fractal dimension of bivariate functions which are of bounded variation in the Arzel\'{a} sense.\\ Fractional calculus is an older topic, dating back almost 300 years which deals with the theory of fractional (non-integer) order derivatives and integrals and their diverse applications. For the relevant literature reader may refer to \cite{Samko}. The fractional calculus of FIF for the case of univariate discussed by Prasad \cite{P} and Gowrisankar \cite{G} which represents a connection between fractional calculus and fractals. But, to the best of our knowledge, there is no research reported in the case  of bivariate FIFs and fractional calculus. In this paper, we established some new results on bivariate FISs and for this, we adopt the construction given by Huo-Jun Ruan and Qiang Xu\cite{Ruan}. By adopting the construction gievn by \cite{Ruan}, we deduce a different kind of bivariate FIFs and FISs and throughout the article, we deal with this.\\
The organization of the article from section 3 to section 7, is as follows: In Section 3, we give construction of bivariate FIFs and FISs. The Partial integral of FIFs and Partial differentiable FISs is given in section 4. In Section 5, we give the definitions of mixed Riemann-Liouville fractional integral and derivative on a closed rectangle $[a,b] \times [c,d]$ and we show that mixed fractional Riemann-Liouville integral and derivative of bivariate FIFs are again bivariate FIFs corresponding to some IFS. Section 6 and 7 are devoted to the integral transforms of the bivariate FIFs.

\section{Background and preliminaries}
In this section, we give terminology, background information, which relates to this report. For  the details, the reader may refer\cite{Ruan}.\\
Let $I=[a,b]$ and $J=[c,d].$ Given interpolation data $\{(x_i,y_j,z_{ij}) \in \mathbb{R}^3 : i=0,1,\dots,N; j=0,1,\dots,M\}$ such that $a=x_0< x_1 < \dots <x_N=b$ and $c=y_0<y_1<\dots <y_M=d.$ For convenience, we write $\Sigma=\{1,2,\dots,N\},$ $ \Sigma_{N,0}=\{0,1,\dots N \},$ $\partial \Sigma_{N,0}=\{0,N\} $ and int$\Sigma_{N,0}=\{1,2,\dots,N-1\}.$ Similarly, we can define $\Sigma_M, \Sigma_{M,0}, \partial \Sigma_{M,0}$ and int$\Sigma_{M,0}.$ Denote $I_i=[x_{i-1},x_i]$ and $J_j=[y_{j-1},y_j]$ for $i \in \Sigma_N$ and $j \in \Sigma_M.$ For any $i \in \Sigma_N,$ let us define a contractive homeomorphism $u_i:I \rightarrow I_i$ satisfying
 \begin{equation}\label{E1}
    \begin{aligned}
    u_i(x_0)=x_{i-1}~ \text{and}~ u_i(x_N)=x_i, ~~\text{when i is odd},
     \end{aligned}
       \end{equation}
       \begin{equation}\label{E2}
          \begin{aligned}
          u_i(x_0)=x_i ~ \text{and}~ u_i(x_N)=x_{i-1},~~ \text{when i is even},
           \end{aligned}
             \end{equation}
    \begin{equation}
             \begin{aligned}
             |u_i(x_1)-u_i(x_2)| \le \alpha_i|x_1 -x_2|, ~~~~ \forall ~~x_1, x_2 \in I,
              \end{aligned}
                \end{equation}

 where $0 < \alpha_i < 1$ is any given constant. One can easily see that,
 \begin{equation*}
  u_1(x_0)=x_0, u_1(x_N)=u_2(x_N)= x_1, u_2(x_0)=u_3(x_0)= x_2~ \text{and so on.}
 \end{equation*}
  Similarly, for any $j \in \Sigma_M,$ let us define a contractive homeomorphism $v_j:J \rightarrow J_j$ satisfying

   \begin{equation}\label{E3}
                 \begin{aligned}
                 v_j(y_0)=y_{j-1}~ \text{and}~ v_j(y_M)=y_j,\text{ when j is odd},
                  \end{aligned}
                    \end{equation}
   \begin{equation}\label{E4}
                 \begin{aligned}
                   v_j(y_0)=y_j ~ \text{and}~ v_j(y_M)=y_{j-1}, ~~ \text{when j is even},
                  \end{aligned}
                    \end{equation}
   \begin{equation}
                 \begin{aligned}
                 |v_j(y_1)-v_j(y_2)| \le \beta_j|y_1 -y_2|, ~~~~ \forall ~~y_1, y_2 \in J,
                  \end{aligned}
                    \end{equation}\\
 where $0 < \beta_j < 1$ is any given constant. From the definitions of $u_i$ and $v_j,$ we can easily check that
 \begin{equation*}
                 \begin{aligned}
                  u^{-1}_i(x_i)=u^{-1}_{i+1}(x_i), ~~~~\forall ~~i \in \text{int}\Sigma_{N,0},
                  \end{aligned}
                    \end{equation*}
                    and
  \begin{equation*}
                  \begin{aligned}
                v^{-1}_j(y_j)=v^{-1}_{j+1}(y_j),~~~~ \forall ~~j \in \text{int}\Sigma_{M,0}.
                   \end{aligned}
                     \end{equation*}
 Let us define $\rho: \mathbb{Z}\times \{0,N,M\} \rightarrow \mathbb{Z}$ by
 \[
 \begin{cases}

                    \rho(i,0)= i-1, \rho(i,N)=\rho(i,M)=i, &  \text{ when i is odd }
                \\

                   \rho(i,0)= i,~~ \rho(i,N)= \rho(i,M)=i-1, & \text{ when i is even}.
                   \end{cases}
                     \]

 Denote $K=I \times J \times \mathbb{R}.$ For each $(i,j) \in \Sigma_N \times \Sigma_M,$ $F_{ij}:K \rightarrow \mathbb{R}$ be a continuous function satisfying $ F_{ij}(x_k,y_l,z_{kl})=z_{\rho(i,k),\rho(j,l)},  ~~~~~~ \forall (k,l) \in \partial \Sigma_{N,0} \times  \partial \Sigma_{M,0}$ and  $|F_{ij}(x,y,z')- F_{ij}(x,y,z'')| \le \gamma_{ij} |z' - z''|, ~~~~~~ \forall (x,y) \in I \times J$ and $z',z'' \in \mathbb{R},$ where $0 < \gamma_{ij} < 1$ is any given constant.\\
 Now, for each $(i,j) \in \Sigma_N \times \Sigma_M,$ we define $W_{ij}:K \rightarrow I_i \times J_j \times \mathbb{R}$ as $$ W_{ij}(x,y,z)=\big(u_i(x),v_j(y),F_{ij}(x,y,z)\big).$$
 Then $\{K, W_{ij}: (i,j) \in\Sigma_N \times \Sigma_M \}$ is an IFS.

 \begin{theorem} \cite{Ruan} \label{BSPTHM1}
 Let $\{K, W_{ij}: (i,j) \in\Sigma_N \times \Sigma_M \}$ be the IFS defined as above. Suppose that $\{F_{ij}:(i,j) \in \Sigma_N \times \Sigma_M \}$ satisfies the below matching conditions:
 \begin{enumerate}
 \item for each $i \in \text{int}\Sigma_{N,0}, j \in \Sigma_M$ and $x^*=u^{-1}_i(x_i)=u^{-1}_{i+1}(x_i),$
  $$F_{ij}(x^*,y,z)=F_{i+1,j}(x^*,y,z),  \forall y \in J, z \in \mathbb{R}, $$
 \item  for each $i \in \Sigma_N, j \in \text{int} \Sigma_{M,0}$ and $y^*=v^{-1}_j(y_j)=v^{-1}_{j+1}(y_j),$ $$F_{ij}(x,y^*,z) = F_{i,j+1}(x,y^*,z),~ \forall x \in I, z \in \mathbb{R}.$$
 \end{enumerate}
  Then there is a unique continuous function $f:I \times J \rightarrow \mathbb{R}  $ satisfying $f(x_i,y_j)=z_{ij}$ for each $(i,j) \in \Sigma_{N,0} \times \Sigma_{M,0}$ and $G= \cup_{(i,j) \in \Sigma_{N} \times \Sigma_{M}} W_{ij}(G),$ where $G=\big\{(x,y,f(x,y)):(x,y) \in I \times J \big\}$ will be the graph of $f.$ We say that $G$ is FIS and $f$ is FIF  corresponding to the IFS $\big\{K, W_{ij}: (i,j) \in\Sigma_N \times \Sigma_M \big\}.$
 \end{theorem}

By Theorem \ref{BSPTHM1}, we have the following result.

\begin{theorem}
Let $\{K, W_{ij}: (i,j) \in\Sigma_N \times \Sigma_M \}$ be the IFS introduced as above.
Then there is a unique continuous function $f:I \times J \rightarrow \mathbb{R}  $ satisfying $f(x_i,y_j)=z_{ij}$ for each $(i,j) \in \Sigma_{N,0} \times \Sigma_{M,0}$ and $G= \cup_{(i,j) \in \Sigma_{N,0} \times \Sigma_{M,0}} W_{ij}(G),$ where $G=\{(x,y,f(x,y)):(x,y) \in I \times J \}$ is the graph of $f.$
\end{theorem}

\section{Construction of Fractal Interpolation Surfaces}
We define $F_{ij}: K \rightarrow \mathbb{R} $ as $ F_{ij}(x,y,z)= \alpha z+ q_{ij}(x,y)$, where $q_{ij}(x_k,y_l)=z_{\rho(i,k),\rho(j,l)}- \alpha z_{kl}$ and $| \alpha |< 1.$ Now, for each $(i,j) \in \Sigma_N \times \Sigma_M,$ we define $W_{ij}:K \rightarrow I_i \times J_j \times \mathbb{R}$ as $$ W_{ij}(x,y,z)=(u_i(x),v_j(y),F_{ij}(x,y,z)).$$
One can see that $F_{ij}$ satisfies matching conditions of Theorem \ref{BSPTHM1}.
The set of all continuous functions on $I\times J$ is denoted by $C(I \times J)$. Define Read-Bajraktarevi\'{c} (RB) operator $T: C(I \times J) \rightarrow C(I \times J)$ by $$ (Tf)(x,y)=F_{ij}\Big(u_i^{-1}(x),v_j^{-1}(y),f\big(u_i^{-1}(x),v_j^{-1}(y)\big)\Big), ~~~\forall ~~ (x,y) \in I_i \times J_j, $$ for all $(i,j) \in \Sigma_N \times \Sigma_M.$
The bivariate FIF $f$ in the above definition will be the unique fixed point of $T$. Consequently, $f$ satisfies the self-referential equation:
$$ f(x,y)=F_{ij}\Big(u_i^{-1}(x),v_j^{-1}(y),f\big(u_i^{-1}(x),v_j^{-1}(y)\big)\Big), ~~~\forall ~~ (x,y) \in I_i \times J_j, $$ for all $(i,j) \in \Sigma_N \times \Sigma_M.$ So, here we have
 \begin{equation}\label{Fnleq1}
   f(x,y)= \alpha~ f\big(u_i^{-1}(x),v_j^{-1}(y)\big)+ q_{ij}\big(u_i^{-1}(x),v_j^{-1}(y)\big),
\end{equation}
for $(x,y) \in I_i \times J_j,~ (i,j) \in \Sigma_N \times \Sigma_M$,
where
\begin{equation}\label{Fn}
u_i(x)=a_ix+b_i ~~\text{and}~~ v_j(y)=c_jy+d_j.
\end{equation}
From \ref{E1}, \ref{E2}, \ref{E3},  \ref{E4} and \ref{Fn}, we yield
\[
 \begin{cases}

                    a_i=\frac{x_i-x_{i-1}}{x_N-x_0} ,~ b_i=\frac{x_{i-1}x_N-x_ix_0}{x_N-x_0}, &  \text{ when i is odd }
                \\

                   a_i=\frac{x_{i-1}-x_{i}}{x_N-x_0} ,~ b_i=\frac{x_{i}x_N-x_{i-1}x_0}{x_N-x_0}, & \text{ when i is even}.
                   \end{cases}
                     \]
                     \[
 \begin{cases}

                    c_j=\frac{y_j-y_{j-1}}{y_M-y_0} ,~ d_j=\frac{y_{j-1}y_M-y_jy_0}{y_M-y_0}, &  \text{ when j is odd }
                \\

                   c_j=\frac{y_{j-1}-y_{j}}{y_M-y_0} ,~ d_j=\frac{y_{j}y_M-y_{j-1}y_0}{y_M-y_0}, & \text{ when j is even}.
                   \end{cases}
                     \]
In this article, we deal with the above equations \ref{Fnleq1} and \ref{Fn}.\\

Let $N=M$. We define a bilinear function $q$ on $I\times J$ which satisfy
\begin{equation*}
q(x_k,y_l)=z_{kl}, ~~~\forall (k,l)\in\partial \Sigma_{N,0}\times\partial \Sigma_{M,0}.
\end{equation*}
That is,
\begin{equation*}
        \begin{aligned}
                  q(x,y)=& \frac{1}{\Delta}[(b-x)(d-y)z_{0,0}+(x-a)(d-y)z_{N,0}\\& +(b-x)(y-c)z_{0,M}+(x-a)(y-c)z_{M,N}].
 \end{aligned}
   \end{equation*}
   Where, $\Delta=(b-a)(d-c)$.
 \subsubsection*{Example:} Let $g(x,y)=\text{sin}(x^2+y^2)$ on $[0,1]\times[0,1]$ ; $x,y\in[0,1]$. Let $N=M=4$. Let $x_i=i/N$ and $y_j=j/M$ for all $i\in \Sigma_{N,0}$ and$j\in \Sigma_{M,0}$.
 So interpolation data $\{(x_i,y_j,g(x_i,y_j)) \in \mathbb{R}^3 : i=0,1,2,3,4; j=0,1,2,3,4\}$ such that $0=x_0< x_1 < x_2 < x_3 < x_4=1$ and $0=y_0<y_1< y_2 < y_3 <y_4=1.$ Let $Z= g(x_i,y_j) =(z_{ij})_{(i,j)\in \Sigma_{N,0} \times \Sigma_{M,0}}.$ Then,
 \[
Z=
  \begin{pmatrix}
    0  &  0.0625   & 0.2474  &  0.5333  &  0.8415 \\
     0.0625 &  0.1247  &  0.3074 &   0.5851  &  0.8736\\
      0.2474 &   0.3074   & 0.4794  &  0.7260 &   0.9490\\
       0.5333  &  0.5851  &  0.7260  &  0.9023  &  1.0000\\
        0.8415   & 0.8736  &  0.9490  &  1.0000 &   0.9093
  \end{pmatrix}.
\]
The corresponding FIS is shown in figure 1.

\begin{figure}
\centering
\includegraphics[scale=0.5]{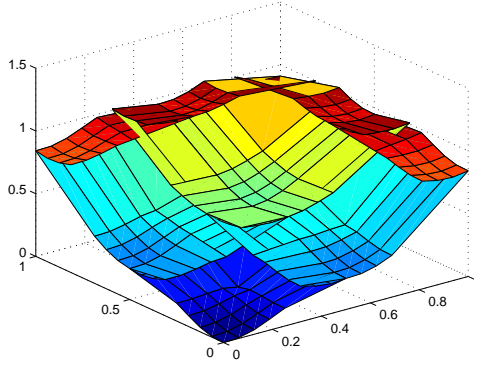}
\caption{ Fractal interpolation surface ($\alpha = 0.2$)}
\end{figure}

\section{Partial Integral of FIS and partial differentiable FIS}

Suppose that $f$ is a FIF defined as above. We show its integral with respect to $x$ is also a fractal interpolation function corresponding to some IFS.
\begin{theorem}
   Let $f$ be the FIF determined by the IFS defined above, then $I_xf$ is FIF assciated with $\{K, \hat{W_{ij}}: (i,j) \in\Sigma_N \times \Sigma_M \},$ where $$ \hat{W}_{ij}(x,y,z)=(u_i(x),v_j(y),\hat{F}_{ij}(x,y,z))$$ and $ \hat{F}_{ij}: K \rightarrow \mathbb{R} $ is defined as $ \hat{F}_{ij}(x,y,z)= \alpha a_i z+ \hat{q}_{ij}(x,y)$ and
   \begin{equation*}
        \begin{aligned}
                    \hat{q}_{ij}(x,y)=&a_i I_x q_{ij}(x,y)+ \int_a ^{u_i(a)} f(s,v_j(y))ds.
       \end{aligned}
   \end{equation*}

   \end{theorem}
   \textbf{Proof.}
Using fixed point equation corresponding to RB operator, we get $f(u_i(x),v_j(y))= \alpha f(x,y) + q_{ij}(x,y), ~ ~~ \forall (x,y) \in I \times J.$
\begin{equation*}
    \begin{aligned}
      I_xf(u_i(x),v_j(y))=& \int_a ^{u_i(x)} f(s,v_j(y))ds\\
      =& \int_a ^{u_i(a)} f(s,v_j(y))ds + \int_{u_i(a)} ^{u_i(x)} f(s,v_j(y))ds\\
      =& \int_a ^{u_i(a)} f(s,v_j(y))ds + a_i \int_{a} ^{x} f(u_i(s),v_j(y))ds\\
       =& \int_a ^{u_i(a)} f(s,v_j(y))ds + a_i \int_{a} ^{x} [\alpha f(s,y) + q_{ij}(s,y)]ds\\
       =& \int_a ^{u_i(a)} f(s,v_j(y))ds + \alpha a_i \int_{a} ^{x}  f(s,y)ds + a_i \int_{a} ^{x} q_{ij}(s,y)ds\\
       =& \alpha a_i I_x f(x,y)+a_i I_x q_{ij}(x,y)+ \int_a ^{u_i(a)} f(s,v_j(y))ds.\\
       =& \alpha a_i I_x f(x,y)+ \hat{q}_{ij}(x,y).
      \end{aligned}
\end{equation*}
This complete the proof.

\begin{remark}
Let $f$ be a differential FIF determined by the IFS defined above. Then, $ \hat{f}=f_x$ if and only if $\hat{f}$ is the FIF associated with $\{K,(u_i(x),v_j(y),\hat{F}_{ij}(x,y,z)) : (i,j) \in\Sigma_N \times \Sigma_M \},$ where $$ \hat{F}_{ij}(x,y,z)= \hat{\alpha} z+ \hat{q}_{ij}(x,y)$$ and $$ \frac{\hat{\alpha}}{a_i}= \alpha, ~~ (\hat{q}_{ij})_x = a_i q_{ij}.$$
\end{remark}

\begin{remark}
In similar way, its integral and derivative with respect to y are also fractal interpolation function corresponding to some IFS.
\end{remark}

\section{Fractional Integral of FIS and Fractional Differentiable FIS}
Mixed Riemann-Liouville fractional order integral and derivative are defined as follows:
\begin{definition}\cite{Samko}
   Let $f$ be a function define on a closed rectangle $[a,b] \times [c,d].$ Assuming that the following integral exists, mixed Riemann-Liouville fractional integral of $f$ is defined by $$ I^{\gamma}f(x,y)=\frac{1}{\Gamma p . \Gamma q} \int_a ^x \int_c ^y (x-s)^{p-1} (y-t)^{q-1}f(s,t)dsdt ,$$ where $\gamma = ( p, q )$ with $ p >0 , q >0.$
   \end{definition}
   \begin{definition}\cite{Samko}
   Let $f$ be a function define on a closed rectangle $[a,b] \times [c,d].$    Assuming that the following integral exists, mixed Riemann-Liouville derivative of $f$ is defined by $$ D^{\gamma}f(x,y)=\frac{1}{\Gamma (1-p) . \Gamma (1-q)} \frac{\partial^2}{\partial{x}\partial{y}} \int_a ^x \int_c ^y (x-s)^{-p} (y-t)^{-q}f(s,t)dsdt ,$$ where $\gamma = ( p, q )$ with $ p >0 , q >0.$
   \end{definition}
We shall prove that fractional integral and derivative of bivariate FIFs are again bivariate FIFs corresponding to some IFS.
   \begin{theorem}
   Let $f$ be the FIF determined by the IFS defined above, then $I^{\gamma}f$ is FIF assciated with $\{K, \hat{W_{ij}}: (i,j) \in\Sigma_N \times \Sigma_M \},$ where $$ \hat{W}_{ij}(x,y,z)=(u_i(x),v_j(y),\hat{F}_{ij}(x,y,z))$$ and $\hat{F}_{ij}: K \rightarrow \mathbb{R} $ is defined as $ \hat{F}_{ij}(x,y,z)= \alpha a_i^p c_j^q z+ \hat{q}_{ij}(x,y)$ where
   \begin{equation*}
        \begin{aligned}
                    \hat{q}_{ij}(x,y)=& \frac{1}{\Gamma p . \Gamma q} \int_a ^{u_i(x)} \int_c ^{v_j(c)} (u_i(x)-s)^{p-1} (v_j(y)-t)^{q-1}f(s,t)dsdt \\
                    &+\frac{c_j^q}{\Gamma p . \Gamma q} \int_a ^{u_i(a)} \int_c ^y (u_i(x)-s)^{p-1} (y-t)^{q-1}f(s,v_j(t))dsdt \\
                    &+\frac{a_i^p c_j^q}{\Gamma p . \Gamma q} \int_a ^x \int_c ^y (x-s)^{p-1} (y-t)^{q-1}q_{ij}(s,t)dsdt .
       \end{aligned}
   \end{equation*}

   \end{theorem}
   \textbf{Proof.}
   Fixed point equation yields $f(u_i(x),v_j(y))= \alpha f(x,y) + q_{ij}(x,y), ~ ~~ \forall (x,y) \in I \times J.$ We have
   \begin{equation*}
    \begin{aligned}
             I^{\gamma}f(u_i(x),v_j(y))=&\frac{1}{\Gamma p . \Gamma q} \int_a ^{u_i(x)} \int_c ^{v_j(y)} (u_i(x)-s)^{p-1} (v_j(y)-t)^{q-1}f(s,t)dsdt\\
              =&\frac{1}{\Gamma p . \Gamma q} \int_a ^{u_i(x)} \int_c ^{v_j(c)} (u_i(x)-s)^{p-1} (v_j(y)-t)^{q-1}f(s,t)dsdt\\ &+\frac{1}{\Gamma p . \Gamma q} \int_a ^{u_i(x)} \int_{v_j(c)} ^{v_j(y)} (u_i(x)-s)^{p-1} (v_j(y)-t)^{q-1}f(s,t)dsdt\\
              =&\frac{1}{\Gamma p . \Gamma q} \int_a ^{u_i(x)} \int_c ^{v_j(c)} (u_i(x)-s)^{p-1} (v_j(y)-t)^{q-1}f(s,t)dsdt\\ &+\frac{c_j}{\Gamma p . \Gamma q} \int_a ^{u_i(x)} \int_{c} ^{y} (u_i(x)-s)^{p-1} (v_j(y)-v_j(w))^{q-1}f(s,v_j(w))dsdw\\
              =&\frac{1}{\Gamma p . \Gamma q} \int_a ^{u_i(x)} \int_c ^{v_j(c)} (u_i(x)-s)^{p-1} (v_j(y)-t)^{q-1}f(s,t)dsdt\\ &+\frac{c_j^q}{\Gamma p . \Gamma q} \int_a ^{u_i(x)} \int_{c} ^{y} (u_i(x)-s)^{p-1} (y-t)^{q-1}f(s,v_j(t))dsdt\\
                      \end{aligned}
   \end{equation*}
              \begin{equation*}
    \begin{aligned}              
              =&\frac{1}{\Gamma p . \Gamma q} \int_a ^{u_i(x)} \int_c ^{v_j(c)} (u_i(x)-s)^{p-1} (v_j(y)-t)^{q-1}f(s,t)dsdt\\ &+\frac{c_j^q}{\Gamma p . \Gamma q} \int_a ^{u_i(a)} \int_{c} ^{y} (u_i(x)-s)^{p-1} (y-t)^{q-1}f(s,v_j(t))dsdt\\ &+ \frac{c_j^q}{\Gamma p . \Gamma q} \int_{u_i(a)} ^{u_i(x)} \int_{c} ^{y} (u_i(x)-s)^{p-1} (y-t)^{q-1}f(s,v_j(t))dsdt\\
              =&\frac{1}{\Gamma p . \Gamma q} \int_a ^{u_i(x)} \int_c ^{v_j(c)} (u_i(x)-s)^{p-1} (v_j(y)-t)^{q-1}f(s,t)dsdt\\ &+\frac{c_j^q}{\Gamma p . \Gamma q} \int_a ^{u_i(a)} \int_{c} ^{y} (u_i(x)-s)^{p-1} (y-t)^{q-1}f(s,v_j(t))dsdt\\ &+ \frac{a_i^p c_j^q}{\Gamma p . \Gamma q} \int_{a} ^{x} \int_{c} ^{y} ( x - s)^{p-1} (y-t)^{q-1}f(u_i(s),v_j(t))dsdt\\
              =&\frac{1}{\Gamma p . \Gamma q} \int_a ^{u_i(x)} \int_c ^{v_j(c)} (u_i(x)-s)^{p-1} (v_j(y)-t)^{q-1}f(s,t)dsdt\\ &+\frac{c_j^q}{\Gamma p . \Gamma q} \int_a ^{u_i(a)} \int_{c} ^{y} (u_i(x)-s)^{p-1} (y-t)^{q-1}f(s,v_j(t))dsdt\\ &+ \frac{\alpha a_i^p c_j^q}{\Gamma p . \Gamma q} \int_{a} ^{x} \int_{c} ^{y} ( x - s)^{p-1} (y-t)^{q-1}f(s,t)dsdt\\ &+ \frac{a_i^p c_j^q}{\Gamma p . \Gamma q} \int_{a} ^{x} \int_{c} ^{y} ( x - s)^{p-1} (y-t)^{q-1} q_{ij}(s,t)dsdt\\
              =&\frac{1}{\Gamma p . \Gamma q} \int_a ^{u_i(x)} \int_c ^{v_j(c)} (u_i(x)-s)^{p-1} (v_j(y)-t)^{q-1}f(s,t)dsdt\\ &+\frac{c_j^q}{\Gamma p . \Gamma q} \int_a ^{u_i(a)} \int_{c} ^{y} (u_i(x)-s)^{p-1} (y-t)^{q-1}f(s,v_j(t))dsdt\\ &+\alpha a_i^p c_j^q ~~I^{\gamma}f(x,y)+ a_i^p c_j^q ~~I^{\gamma}q_{ij}(x,y)\\
              =& \alpha a_i^p c_j^q~~ I^{\gamma}f(x,y)+ \hat{q}_{ij}(x,y).
               \end{aligned}
   \end{equation*}
 Hence the required result follows.

     \begin{theorem}
   Let $f$ be a fractional differentiable FIF determined by the IFS defined above, then $D^{\gamma}f$ is FIF assciated with $\{K, \hat{W_{ij}}: (i,j) \in\Sigma_N \times \Sigma_M \},$ where $$ \hat{W}_{ij}(x,y,z)=(u_i(x),v_j(y),\hat{F}_{ij}(x,y,z))$$ and $\hat{F}_{ij}: K \rightarrow \mathbb{R} $ is defined as $ \hat{F}_{ij}(x,y,z)= \alpha a_i^{(1-p)} c_j^{(1-q)} z+ \hat{q}_{ij}(x,y)$ where
   \begin{equation*}
        \begin{aligned}
                    \hat{q}_{ij}(x,y)=& \frac{1}{\Gamma (1-p) . \Gamma (1-q)}\frac{\partial^2}{\partial{x}\partial{y}} \int_a ^{u_i(x)} \int_c ^{v_j(c)} (u_i(x)-s)^{-p} (v_j(y)-t)^{-q}f(s,t)dsdt \\
                    &+\frac{c_j^{(1-q)}}{\Gamma (1-p) . \Gamma (1-q)}\frac{\partial^2}{\partial{x}\partial{y}} \int_a ^{u_i(a)} \int_{c} ^{y} (u_i(x)-s)^{-p} (y-t)^{-q}f(s,v_j(t))dsdt \\
                    &+\frac{a_i^{(1-p)} c_j^{(1-q)}}{\Gamma (1-p) . \Gamma (1-q)}\frac{\partial^2}{\partial{x}\partial{y}} \int_{a} ^{x} \int_{c} ^{y} ( x - s)^{-p} (y-t)^{-q} q_{ij}(s,t)dsdt .
       \end{aligned}
   \end{equation*}

   \end{theorem}
   \textbf{Proof.}
   Fixed point equation yields $f(u_i(x),v_j(y))= \alpha f(x,y) + q_{ij}(x,y), ~ ~~ \forall (x,y) \in I \times J.$ We have
   \begin{equation*}
    \begin{aligned}
             D^{\gamma}f(u_i(x),v_j(y))=&\frac{1}{\Gamma (1-p) . \Gamma (1-q)} \frac{\partial^2}{\partial{x}\partial{y}} \int_a ^{u_i(x)} \int_c ^{v_j(y)} (u_i(x)-s)^{-p} (v_j(y)-t)^{-q}f(s,t)dsdt\\
              =&\frac{1}{\Gamma (1-p) . \Gamma (1-q)}\frac{\partial^2}{\partial{x}\partial{y}} \big[\int_a ^{u_i(x)} \int_c ^{v_j(c)} (u_i(x)-s)^{-p} (v_j(y)-t)^{-q}f(s,t)dsdt\\ &+ \int_a ^{u_i(x)} \int_{v_j(c)} ^{v_j(y)} (u_i(x)-s)^{-p} (v_j(y)-t)^{-q}f(s,t)dsdt \big]\\
              =&\frac{1}{\Gamma (1-p) . \Gamma (1-q)}\frac{\partial^2}{\partial{x}\partial{y}} \big[ \int_a ^{u_i(x)} \int_c ^{v_j(c)} (u_i(x)-s)^{-p} (v_j(y)-t)^{-q}f(s,t)dsdt\\ &+{c_j} \int_a ^{u_i(x)} \int_{c} ^{y} (u_i(x)-s)^{-p} (v_j(y)-v_j(w))^{-q}f(s,v_j(w))dsdw \big] \\                   
              =&\frac{1}{\Gamma (1-p) . \Gamma (1-q)}\frac{\partial^2}{\partial{x}\partial{y}} \big[ \int_a ^{u_i(x)} \int_c ^{v_j(c)} (u_i(x)-s)^{-p} (v_j(y)-t)^{-q}f(s,t)dsdt\\ &+{c_j^{(1-q)}} \int_a ^{u_i(x)} \int_{c} ^{y} (u_i(x)-s)^{-p} (y-t)^{-q}f(s,v_j(t))dsdt \big]\\
              =&\frac{1}{\Gamma (1-p) . \Gamma (1-q)}\frac{\partial^2}{\partial{x}\partial{y}} \big[ \int_a ^{u_i(x)} \int_c ^{v_j(c)} (u_i(x)-s)^{-p} (v_j(y)-t)^{-q}f(s,t)dsdt\\ &+{c_j^{(1-q)}} \int_a ^{u_i(a)} \int_{c} ^{y} (u_i(x)-s)^{-p} (y-t)^{-q}f(s,v_j(t))dsdt\\ &+ {c_j^{(1-q)}} \int_{u_i(a)} ^{u_i(x)} \int_{c} ^{y} (u_i(x)-s)^{-p} (y-t)^{-q}f(s,v_j(t))dsdt \big]\\
              =&\frac{1}{\Gamma (1-p) . \Gamma (1-q)}\frac{\partial^2}{\partial{x}\partial{y}} \big[ \int_a ^{u_i(x)} \int_c ^{v_j(c)} (u_i(x)-s)^{-p} (v_j(y)-t)^{-q}f(s,t)dsdt\\ &+{c_j^{(1-q)}} \int_a ^{u_i(a)} \int_{c} ^{y} (u_i(x)-s)^{-p} (y-t)^{-q}f(s,v_j(t))dsdt\\ &+ {a_i^{(1-p)} c_j^{(1-q)}} \int_{a} ^{x} \int_{c} ^{y} ( x - s)^{-p} (y-t)^{-q}f(u_i(s),v_j(t))dsdt \big]\\
              =&\frac{1}{\Gamma (1-p) . \Gamma (1-q)}\frac{\partial^2}{\partial{x}\partial{y}} \big[ \int_a ^{u_i(x)} \int_c ^{v_j(c)} (u_i(x)-s)^{-p} (v_j(y)-t)^{-q}f(s,t)dsdt\\ &+{c_j^{(1-q)}} \int_a ^{u_i(a)} \int_{c} ^{y} (u_i(x)-s)^{-p} (y-t)^{-q}f(s,v_j(t))dsdt\\ &+ {\alpha a_i^{(1-p)} c_j^{(1-q)}} \int_{a} ^{x} \int_{c} ^{y} ( x - s)^{-p} (y-t)^{-q}f(s,t)dsdt\\ &+ {a_i^{(1-p)} c_j^{(1-q)}} \int_{a} ^{x} \int_{c} ^{y} ( x - s)^{-p} (y-t)^{-q} q_{ij}(s,t)dsdt \big]\\
              =&\frac{1}{\Gamma (1-p) . \Gamma (1-q)}\frac{\partial^2}{\partial{x}\partial{y}} \int_a ^{u_i(x)} \int_c ^{v_j(c)} (u_i(x)-s)^{-p} (v_j(y)-t)^{-q}f(s,t)dsdt\\ &+ \frac{c_j^{(1-q)}}{\Gamma (1-p) . \Gamma (1-q)}\frac{\partial^2}{\partial{x}\partial{y}} \int_a ^{u_i(a)} \int_{c} ^{y} (u_i(x)-s)^{-p} (y-t)^{-q}f(s,v_j(t))dsdt\\ &+\alpha a_i^{(1-p)} c_j^{(1-q)} ~~D^{\gamma}f(x,y)+ a_i^{(1-p)} c_j^{(1-q)} ~~D^{\gamma}q_{ij}(x,y)\\
              =& \alpha a_i^{(1-p)} c_j^{(1-q)}~~ D^{\gamma}f(x,y)+ \hat{q}_{ij}(x,y).
               \end{aligned}
   \end{equation*}
From which the required result follows.

 \section{Integral Transforms of bivariate FIFs}
 In this section and in the next section, we present integral transforms corresponding to bivariate FIFs and also remark (see remark \ref{R}) that these integral transforms again form bivariate FIFs.\\
  The general integral transform of a "well-behaved" function $g$ is defined as $$ \hat g(s,t)= \int_{ \mathbb{R}^2} K(x,y,s,t)g(x,y)dxdy ,$$ where $K(x,y,s,t)$ is a suitable function, referred to as kernel of the transformation. Here, we have
  \begin{equation*}
             \begin{aligned}
            \hat f(s,t)=& \int_{ \mathbb{R}^2} K(x,y,s,t)f(x,y)dxdy\\
       =& \int_{ I \times J} K(x,y,s,t)f(x,y)dxdy\\
             =& \Sigma_{j=1}^{M} \Sigma_{i=1}^{N} \int_{ I_i \times J_j} K(x,y,s,t)[q_{ij}(u_i^{-1}(x),v_j^{-1}(y))+ \alpha f(u_i^{-1}(x),v_j^{-1}(y))] dxdy\\
              =& \Sigma_{j=1}^{M} \Sigma_{i=1}^{N} \int_{ I_i \times J_j} K(x,y,s,t)q_{ij}(u_i^{-1}(x),v_j^{-1}(y)) dxdy \\&+  \Sigma_{j=1}^{M} \Sigma_{i=1}^{N} \int_{ I_i \times J_j}\alpha K(x,y,s,t) f(u_i^{-1}(x),v_j^{-1}(y))  dxdy.
                  \end{aligned}
                 \end{equation*}
  With the change of variable $u_i^{-1}(x)=z$ and $v_j^{-1}(y)=w,$ we have
  \begin{equation*}
              \begin{aligned}
              \hat f(s,t)
              =&\Sigma_{j=1}^{M} \Sigma_{i=1}^{N} \int_{ I_i \times J_j} a_i c_j K(u_i(z),v_j(w),s,t) q_{ij}(z,w) dzdw\\&+ \Sigma_{j=1}^{M} \Sigma_{i=1}^{N} \int_{ I_i \times J_j}  \alpha a_i c_j K(u_i(z),v_j(w),s,t) f(z,w) dzdw\\
              =&\Sigma_{j=1}^{M} \Sigma_{i=1}^{N} \int_{ I_i \times J_j}   a_i c_j K(u_i(x),v_j(y),s,t) q_{ij}(x,y) dxdy\\&+ \Sigma_{j=1}^{M} \Sigma_{i=1}^{N} \int_{ I_i \times J_j}  \alpha a_i c_j K(u_i(x),v_j(y),s,t) f(x,y) dxdy\\
               =&\Sigma_{j=1}^{M} \Sigma_{i=1}^{N}a_i c_j \int_{ I_i \times J_j} K(u_i(x),v_j(y),s,t) q_{ij}(x,y) dxdy\\&+ \Sigma_{j=1}^{M} \Sigma_{i=1}^{N} \alpha a_i c_j \int_{ I_i \times J_j}   K(u_i(x),v_j(y),s,t) f(x,y) dxdy.
                  \end{aligned}
                  \end{equation*}
 We study the integral transforms of fractal functions with some special choice of kernel functions.\\
 (1) Laplace Transform:\\
 The two-dimensional Laplace transform of $g$ is defined as $$ \mathcal{L}(g)(s,t)=\int_{0}^{\infty} \int_{0}^{\infty} g(x,y) \exp(-xs-yt)dxdy.  $$
 Here $K(x,y,s,t)= \exp(-xs-yt)$ for $x ,y ,s,t>0.$ We choose $I=[a,b]$ and $J=[c,d]$ such that both $a,c >0.$ Then
 \begin{equation*}
              \begin{aligned}
              \hat f(s,t)
              =& \Sigma_{j=1}^{M} \Sigma_{i=1}^{N}a_i c_j \int_{ I_i \times J_j} K(u_i(x),v_j(y),s,t) q_{ij}(x,y) dxdy\\&+ \Sigma_{j=1}^{M} \Sigma_{i=1}^{N} \alpha a_i c_j \int_{ I_i \times J_j}   K(u_i(x),v_j(y),s,t) f(x,y) dxdy\\
             =& \Sigma_{j=1}^{M} \Sigma_{i=1}^{N}a_i c_j \int_{ I_i \times J_j} \exp(-(u_i(x)s+v_j(y)t) q_{ij}(x,y) dxdy\\&+ \Sigma_{j=1}^{M} \Sigma_{i=1}^{N}  \alpha a_i c_j  \int_{ I_i \times J_j} \exp(-(u_i(x)s+v_j(y)t)) f(x,y) dxdy\\
             =& \Sigma_{j=1}^{M} \Sigma_{i=1}^{N}a_i c_j \int_{ I_i \times J_j} \exp(-((a_ix+b_i)s+(c_jy+d_j)t)) q_{ij}(x,y) dxdy\\&+ \Sigma_{j=1}^{M} \Sigma_{i=1}^{N}  \alpha a_i c_j  \int_{ I_i \times J_j} \exp(-((a_ix+b_i)s+(c_jy+d_j)t)) f(x,y) dxdy\\
             =&\Sigma_{j=1}^{M} \Sigma_{i=1}^{N}  a_i c_j \exp (-(b_is+d_jt))\int_{ I_i \times J_j}   \exp(-(xa_is+yc_jt)) q_{ij}(x,y) dxdy \\&+ \Sigma_{j=1}^{M} \Sigma_{i=1}^{N} \alpha a_i c_j \exp (-(b_is+d_jt))\int_{ I_i \times J_j}   \exp(-(xa_is+yc_jt)) f(x,y) dxdy\\
              =& \Sigma_{j=1}^{M} \Sigma_{i=1}^{N}  a_i c_j \exp (-(b_is+d_jt))q_{ij}(a_is,c_jt)\\&+ \Sigma_{j=1}^{M} \Sigma_{i=1}^{N} \alpha a_i c_j \exp (-(b_is+d_jt))f(a_is,c_jt)\\
              =& \hat q_{ij}(s,t)+ \Sigma_{j=1}^{M} \Sigma_{i=1}^{N} \alpha a_i c_j \exp (-(b_is+d_jt))f(a_is,c_jt) .\\
                  \end{aligned}
                  \end{equation*}
 Where, $\hat q_{ij}(s,t)=\Sigma_{j=1}^{M} \Sigma_{i=1}^{N}  a_i c_j \exp (-(b_is+d_jt))q_{ij}(a_is,c_jt).$ \\
 
 (2) Laplace Carson Transform:\\
 The two-dimensional Laplace carson transform of $g$ is defined as $$ \mathcal{L}(g)(s,t)=st\int_{0}^{\infty} \int_{0}^{\infty} g(x,y) \exp(-xs-yt)dxdy.  $$
 Here $K(x,y,s,t)= st\exp(-xs-yt)$ for $x ,y ,s,t>0.$ We choose $I=[a,b]$ and $J=[c,d]$ such that both $a,c >0.$ Then 
 \begin{equation*}
              \begin{aligned} 
              \hat f(s,t)
              =& \Sigma_{j=1}^{M} \Sigma_{i=1}^{N}a_i c_j \int_{ I_i \times J_j} K(u_i(x),v_j(y),s,t) q_{ij}(x,y) dxdy\\&+ \Sigma_{j=1}^{M} \Sigma_{i=1}^{N} \alpha a_i c_j \int_{ I_i \times J_j}   K(u_i(x),v_j(y),s,t) f(x,y) dxdy\\
             =& \Sigma_{j=1}^{M} \Sigma_{i=1}^{N}a_i c_j st\int_{ I_i \times J_j} \exp(-(u_i(x)s+v_j(y)t) q_{ij}(x,y) dxdy\\&+ \Sigma_{j=1}^{M} \Sigma_{i=1}^{N}  \alpha a_i c_j st \int_{ I_i \times J_j} \exp(-(u_i(x)s+v_j(y)t)) f(x,y) dxdy\\
             =& \Sigma_{j=1}^{M} \Sigma_{i=1}^{N}a_i c_jst \int_{ I_i \times J_j} \exp(-((a_ix+b_i)s+(c_jy+d_j)t)) q_{ij}(x,y) dxdy\\&+ \Sigma_{j=1}^{M} \Sigma_{i=1}^{N}  \alpha a_i c_j st \int_{ I_i \times J_j} \exp(-((a_ix+b_i)s+(c_jy+d_j)t)) f(x,y) dxdy\\
           \end{aligned}
                  \end{equation*}
            \begin{equation*}
              \begin{aligned}
             =&\Sigma_{j=1}^{M} \Sigma_{i=1}^{N}  a_i c_j st\exp (-(b_is+d_jt))\int_{ I_i \times J_j}   \exp(-(xa_is+yc_jt)) q_{ij}(x,y) dxdy \\&+ \Sigma_{j=1}^{M} \Sigma_{i=1}^{N} \alpha a_i c_j  st\exp (-(b_is+d_jt))\int_{ I_i \times J_j}   \exp(-(xa_is+yc_jt)) f(x,y) dxdy\\            
              =& \Sigma_{j=1}^{M} \Sigma_{i=1}^{N}  a_i c_jst \exp (-(b_is+d_jt))q_{ij}(a_is,c_jt)\\&+ \Sigma_{j=1}^{M} \Sigma_{i=1}^{N} \alpha a_i c_j st\exp (-(b_is+d_jt))f(a_is,c_jt) \\
              =&\hat q_{ij}(s,t) + \Sigma_{j=1}^{M} \Sigma_{i=1}^{N} \alpha a_i c_j st\exp (-(b_is+d_jt))f(a_is,c_jt) .
                  \end{aligned}
                  \end{equation*}  
Where, $\hat q_{ij}(s,t)=\Sigma_{j=1}^{M} \Sigma_{i=1}^{N}  a_i c_jst \exp (-(b_is+d_jt))q_{ij}(a_is,c_jt).$ \\

 (3) Fourier Transform:\\
 For the kernel $K(x,y,s,t)= \exp(-2 \pi k(xs+yt)),$ where $k = \sqrt{-1}.$ We apply the same process, as above, to evaluate Fourier transform of fractal interpolation function.
 \begin{equation*}
              \begin{aligned}
              \hat f(s,t)
              =& \Sigma_{j=1}^{M} \Sigma_{i=1}^{N}a_i c_j \int_{ I_i \times J_j} K(u_i(x),v_j(y),s,t) q_{ij}(x,y) dxdy\\&+ \Sigma_{j=1}^{M} \Sigma_{i=1}^{N} \alpha a_i c_j \int_{ I_i \times J_j}   K(u_i(x),v_j(y),s,t) f(x,y) dxdy\\
              =& \Sigma_{j=1}^{M} \Sigma_{i=1}^{N}   a_i c_j  \int_{ I_i \times J_j} \exp(-2\pi k(u_i(x)s+v_j(y)t)) q_{ij}(x,y) dxdy\\&+ \Sigma_{j=1}^{M} \Sigma_{i=1}^{N}  \alpha a_i c_j  \int_{ I_i \times J_j} \exp(-2\pi k(u_i(x)s+v_j(y)t)) f(x,y) dxdy\\
              =& \Sigma_{j=1}^{M} \Sigma_{i=1}^{N}   a_i c_j  \int_{ I_i \times J_j} \exp(-2\pi k((a_ix+b_i)s+(c_jy+d_j)t)) q_{ij}(x,y) dxdy\\&+ \Sigma_{j=1}^{M} \Sigma_{i=1}^{N}  \alpha a_i c_j  \int_{ I_i \times J_j} \exp(-2\pi k((a_ix+b_i)s+(c_jy+d_j)t)) f(x,y) dxdy\\
              =& \Sigma_{j=1}^{M} \Sigma_{i=1}^{N}  a_i c_j \exp (-2 \pi k(b_is+d_jt))\int_{ I_i \times J_j}   \exp(-2 \pi k(xa_is+yc_jt)) q_{ij}(x,y) dxdy\\\end{aligned}
                  \end{equation*}
                  \begin{equation*}
              \begin{aligned}
              &+ \Sigma_{j=1}^{M} \Sigma_{i=1}^{N} \alpha a_i c_j \exp (-2 \pi k(b_is+d_jt))\int_{ I_i \times J_j}   \exp(-2 \pi k(xa_is+yc_jt)) f(x,y) dxdy\\
                           =& \Sigma_{j=1}^{M} \Sigma_{i=1}^{N}  a_i c_j \exp( -2 \pi k(b_is+d_jt))q_{ij}(a_is,c_jt)\\&+ \Sigma_{j=1}^{M} \Sigma_{i=1}^{N} \alpha a_i c_j \exp( -2 \pi k(b_is+d_jt))f(a_is,c_jt)\\
                           =&\hat q_{ij}(s,t) + \Sigma_{j=1}^{M} \Sigma_{i=1}^{N} \alpha a_i c_j \exp( -2 \pi k(b_is+d_jt))f(a_is,c_jt) .
                  \end{aligned}
                  \end{equation*}
Where, $\hat q_{ij}(s,t)=\Sigma_{j=1}^{M} \Sigma_{i=1}^{N}  a_i c_j \exp( -2 \pi k(b_is+d_jt))q_{ij}(a_is,c_jt)$.\\
\section{Fractional Order Integral Transforms of bivariate FIFs}
(1) Fractional Double Sumudu Transform:
\begin{definition}\cite{OM} Assuming that $f(x,y)$ vanishes for negative values of $x$ and $y$. The Fractional double Sumudu transform of  $f(x,y)$ of order $\lambda>0$ is defined by
\begin{equation*}
 \mathcal{G}^2_{\lambda}(s,t)=\int_{0}^{\infty} \int_{0}^{\infty} E_{\lambda}(-(x+y)^{\lambda})f(sx,ty)(dx)^\lambda(dy)^\lambda.
\end{equation*}
Where $s,t\in\mathbb{C}$, $ x,y>0$ and $E_{\lambda}(x)=\sum_{m=0}^{\infty}\frac{x^m}{\Gamma(\lambda m +1)}$ is the \textit{Mittag-Leffler function.}
\end{definition}
By applying the same process and terminology, as above and taking  $s,t\neq 0$, we have
\begin{equation*}
             \begin{aligned}
             \hat f(s,t)
              =& \Sigma_{j=1}^{M} \Sigma_{i=1}^{N} \int_{ I_i \times J_j}E_{\lambda}(-(x+y)^{\lambda}) q_{ij}(u_i^{-1}(sx),v_j^{-1}(ty)) (dx)^\lambda (dy)^\lambda \\&+  \Sigma_{j=1}^{M} \Sigma_{i=1}^{N} \int_{ I_i \times J_j}\alpha E_{\lambda}(-(x+y)^{\lambda}) f(u_i^{-1}(sx),v_j^{-1}(ty)) (dx)^\lambda (dy)^\lambda.
                  \end{aligned}
                 \end{equation*}
With the change of variable $u_i^{-1}(sx)=z$ and $v_j^{-1}(ty)=w,$ we have

 \begin{equation*}
             \begin{aligned}
             \hat f(s,t)
              =& \Sigma_{j=1}^{M} \Sigma_{i=1}^{N} \int_{ I_i \times J_j}(\frac{a_i}{s})^\lambda (\frac{c_i}{t})^\lambda E_{\lambda}(-(\frac{u_i(z)}{s}+\frac{v_j(w)}{t})^{\lambda}) q_{ij}(z,w) (dz)^\lambda (dw)^\lambda \\&+  \Sigma_{j=1}^{M} \Sigma_{i=1}^{N} \int_{ I_i \times J_j}\alpha (\frac{a_i}{s})^\lambda (\frac{c_i}{t})^\lambda E_{\lambda}(-(\frac{u_i(z)}{s}+\frac{v_j(w)}{t})^{\lambda}) f(z,w) (dz)^\lambda (dw)^\lambda\\
              =& \Sigma_{j=1}^{M} \Sigma_{i=1}^{N} \int_{ I_i \times J_j}(\frac{a_i}{s})^\lambda (\frac{c_i}{t})^\lambda E_{\lambda}(-(\frac{u_i(x)}{s}+\frac{v_j(y)}{t})^{\lambda}) q_{ij}(x,y) (dx)^\lambda (dy)^\lambda \\&+  \Sigma_{j=1}^{M} \Sigma_{i=1}^{N} \int_{ I_i \times J_j}\alpha (\frac{a_i}{s})^\lambda (\frac{c_i}{t})^\lambda E_{\lambda}(-(\frac{u_i(x)}{s}+\frac{v_j(y)}{t})^{\lambda}) f(x,y) (dx)^\lambda (dy)^\lambda.
                  \end{aligned}
                 \end{equation*}
By applying the Mittag-Leffler property, $E_{\lambda}(k(x+y)^\lambda)=E_{\lambda}(kx^\lambda)E_{\lambda}(ky^\lambda)$, we have
 \begin{equation*}
             \begin{aligned}
             \hat f(s,t)
             =& \Sigma_{j=1}^{M} \Sigma_{i=1}^{N} \int_{ I_i \times J_j}(\frac{a_i}{s})^\lambda (\frac{c_i}{t})^\lambda E_{\lambda}(-(\frac{u_i(x)}{s})^\lambda)E_{\lambda}(-(\frac{v_j(y)}{t})^\lambda) q_{ij}(x,y) (dx)^\lambda (dy)^\lambda \\&+  \Sigma_{j=1}^{M} \Sigma_{i=1}^{N} \int_{ I_i \times J_j}\alpha (\frac{a_i}{s})^\lambda (\frac{c_i}{t})^\lambda E_{\lambda}(-(\frac{u_i(x)}{s})^\lambda)E_{\lambda}(-(\frac{v_j(y)}{t})^\lambda) f(x,y) (dx)^\lambda (dy)^\lambda\\
             =& \Sigma_{j=1}^{M} \Sigma_{i=1}^{N} \int_{ I_i \times J_j}(\frac{a_i}{s})^\lambda (\frac{c_i}{t})^\lambda E_{\lambda}(-(\frac{a_ix+b_i}{s})^\lambda)E_{\lambda}(-(\frac{c_jy+d_j}{t})^\lambda) q_{ij}(x,y) (dx)^\lambda (dy)^\lambda \\&+  \Sigma_{j=1}^{M} \Sigma_{i=1}^{N} \int_{ I_i \times J_j}\alpha (\frac{a_i}{s})^\lambda (\frac{c_i}{t})^\lambda E_{\lambda}(-(\frac{a_ix+b_i}{s})^\lambda)E_{\lambda}(-(\frac{c_jy+d_j)}{t})^\lambda) f(x,y) (dx)^\lambda (dy)^\lambda\\
             =& \Sigma_{j=1}^{M} \Sigma_{i=1}^{N} (\frac{a_i}{s})^\lambda (\frac{c_i}{t})^\lambda  E_{\lambda}(-(\frac{b_i}{s}+\frac{d_j}{t})^\lambda) \int_{ I_i \times J_j} E_{\lambda}(-(\frac{a_ix}{s}+\frac{c_jy}{t})^\lambda) q_{ij}(x,y) (dx)^\lambda (dy)^\lambda \\&+  \Sigma_{j=1}^{M} \Sigma_{i=1}^{N} \alpha (\frac{a_i}{s})^\lambda (\frac{c_i}{t})^\lambda E_{\lambda}(-(\frac{b_i}{s}+\frac{d_j}{t})^\lambda) \int_{ I_i \times J_j} E_{\lambda}(-(\frac{a_ix}{s}+\frac{c_jy}{t})^\lambda) f(x,y) (dx)^\lambda (dy)^\lambda\\
             =& \Sigma_{j=1}^{M} \Sigma_{i=1}^{N} (\frac{a_i}{s})^\lambda (\frac{c_i}{t})^\lambda  E_{\lambda}(-(\frac{b_i}{s}+\frac{d_j}{t})^\lambda) \int_{ I_i \times J_j} E_{\lambda}(-(x+y)^\lambda) q_{ij}(\frac{sx}{a_i},\frac{ty}{c_j}) (dx)^\lambda (dy)^\lambda \\&+  \Sigma_{j=1}^{M} \Sigma_{i=1}^{N} \alpha (\frac{a_i}{s})^\lambda (\frac{c_i}{t})^\lambda E_{\lambda}(-(\frac{b_i}{s}+\frac{d_j}{t})^\lambda) \int_{ I_i \times J_j} E_{\lambda}(-(x+y)^\lambda) f(\frac{sx}{a_i},\frac{ty}{c_j}) (dx)^\lambda (dy)^\lambda\\
             =& \Sigma_{j=1}^{M} \Sigma_{i=1}^{N} (\frac{a_i}{s})^\lambda (\frac{c_i}{t})^\lambda  E_{\lambda}(-(\frac{b_i}{s}+\frac{d_j}{t})^\lambda)  q_{ij}(\frac{s}{a_i},\frac{t}{c_j}) \\&+  \Sigma_{j=1}^{M} \Sigma_{i=1}^{N} \alpha (\frac{a_i}{s})^\lambda (\frac{c_i}{t})^\lambda E_{\lambda}(-(\frac{b_i}{s}+\frac{d_j}{t})^\lambda)) f(\frac{s}{a_i},\frac{t}{c_j})\\
             =& \hat q_{ij}(s,t)+\Sigma_{j=1}^{M} \Sigma_{i=1}^{N} \alpha (\frac{a_i}{s})^\lambda (\frac{c_i}{t})^\lambda E_{\lambda}(-(\frac{b_i}{s}+\frac{d_j}{t})^\lambda)) f(\frac{s}{a_i},\frac{t}{c_j}).
  \end{aligned}
                 \end{equation*}
 Where, $\hat q_{ij}(s,t)=\Sigma_{j=1}^{M} \Sigma_{i=1}^{N} (\frac{a_i}{s})^\lambda (\frac{c_i}{t})^\lambda  E_{\lambda}(-(\frac{b_i}{s}+\frac{d_j}{t})^\lambda)  q_{ij}(\frac{s}{a_i},\frac{t}{c_j}).$\\
 
 (2) Double Laplace Transform of Fractional Order:
 \begin{definition}\cite{OM} If $f(x,y)$ is a function where $x,y>0$, then double Laplace transform of fractional order($\lambda>0$) of function $f(x,y)$ is defined by
 \begin{equation*}
 \mathcal{F}^2_{\lambda}(s,t)=\int_{0}^{\infty} \int_{0}^{\infty} E_{\lambda}(-(sx+ty)^{\lambda})f(x,y)(dx)^\lambda(dy)^\lambda.
\end{equation*}
Where $s,t\in\mathbb{C}$ and $E_{\lambda}(x)$ is the \textit{Mittag-Leffler function.}
\end{definition}
For $\lambda=1$, it turns into double Laplace transform.\\
By applying same process and terminology, as above, we have
\begin{equation*}
              \begin{aligned}
              \hat f(s,t)
               =&\Sigma_{j=1}^{M} \Sigma_{i=1}^{N}a_i^\lambda c_j^\lambda \int_{ I_i \times J_j} K(u_i(x),v_j(y),s,t) q_{ij}(x,y) (dx)^\lambda (dy)^\lambda\\&+ \Sigma_{j=1}^{M} \Sigma_{i=1}^{N} \alpha a_i^\lambda c_j^\lambda \int_{ I_i \times J_j}   K(u_i(x),v_j(y),s,t) f(x,y) (dx)^\lambda (dy)^\lambda.
                  \end{aligned}
                  \end{equation*}
Here $K(x,y,s,t)=E_{\lambda}(-(sx+ty)^{\lambda})$ and by applying the Mittag-Leffler property, we have
\begin{equation*}
              \begin{aligned}
              \hat f(s,t)
               =&\Sigma_{j=1}^{M} \Sigma_{i=1}^{N}a_i^\lambda c_j^\lambda \int_{ I_i \times J_j} E_{\lambda}(-(u_i(x)s)^\lambda)E_{\lambda}(-(v_j(y)t)^{\lambda}) q_{ij}(x,y) (dx)^\lambda (dy)^\lambda\\&+ \Sigma_{j=1}^{M} \Sigma_{i=1}^{N} \alpha a_i^\lambda c_j^\lambda \int_{ I_i \times J_j}   E_{\lambda}(-(u_i(x)s)^\lambda)E_{\lambda}(-(v_j(y)t)^{\lambda}) f(x,y) (dx)^\lambda (dy)^\lambda\\
               =&\Sigma_{j=1}^{M} \Sigma_{i=1}^{N}a_i^\lambda c_j^\lambda \int_{ I_i \times J_j} E_{\lambda}(-((a_ix+b_i)s)^\lambda)E_{\lambda}(-((c_jy+d_j)t)^{\lambda}) q_{ij}(x,y) (dx)^\lambda (dy)^\lambda\\&+ \Sigma_{j=1}^{M} \Sigma_{i=1}^{N} \alpha a_i^\lambda c_j^\lambda \int_{ I_i \times J_j}   E_{\lambda}(-((a_ix+b_i)s)^\lambda)E_{\lambda}(-((c_jy+d_j)t)^{\lambda}) f(x,y) (dx)^\lambda (dy)^\lambda\\
               =&\Sigma_{j=1}^{M} \Sigma_{i=1}^{N}a_i^\lambda c_j^\lambda  E_{\lambda}(-(b_is+d_jt)^\lambda)\int_{ I_i \times J_j} E_{\lambda}(-((a_ixs+c_jyt)^\lambda) q_{ij}(x,y) (dx)^\lambda (dy)^\lambda\\&+ \Sigma_{j=1}^{M} \Sigma_{i=1}^{N} \alpha a_i^\lambda c_j^\lambda E_{\lambda}(-(b_is+d_jt)^\lambda) \int_{ I_i \times J_j} E_{\lambda}(-((a_ixs+c_jyt)^\lambda  f(x,y) (dx)^\lambda (dy)^\lambda\\
               =&\Sigma_{j=1}^{M} \Sigma_{i=1}^{N}a_i^\lambda c_j^\lambda  E_{\lambda}(-(b_is+d_jt)^\lambda)q_{ij}(a_is,c_jt)\\&+ \Sigma_{j=1}^{M} \Sigma_{i=1}^{N} \alpha a_i^\lambda c_j^\lambda E_{\lambda}(-(b_is+d_jt)^\lambda) f(a_is,c_jt)\\
               =&\hat q_{ij}(s,t)+\Sigma_{j=1}^{M} \Sigma_{i=1}^{N} \alpha a_i^\lambda c_j^\lambda E_{\lambda}(-(b_is+d_jt)^\lambda) f(a_is,c_jt).
                  \end{aligned}
                  \end{equation*}
Where, $\hat q_{ij}(s,t)=\Sigma_{j=1}^{M} \Sigma_{i=1}^{N}a_i^\lambda c_j^\lambda  E_{\lambda}(-(b_is+d_jt)^\lambda)q_{ij}(a_is,c_jt)$.\\
\begin{remark}\label{R}
Integral transforms and fractional order integral transforms defined above are also bivariate FIFs for given $\hat q_{ij}(x,y)$ and $|\alpha|<1$ corresponding to some IFS.
\end{remark}
\subsubsection*{Conclusion}
Let us reinforce that the bivariate case is not merely presented as
an extension of the univariate case; but we explored results on fractal functions and fractional calculus. We give several basic and advanced results on bivariate FIFs. In short, we conclude that the various type of transforms of bivariate FIFs are again bivariate FIFs corresponding to some IFS.

\subsection*{Acknowledgements} The work of first author (SC) is financially supported by the CSIR, India with fellowship number 09/1058(0012)/2018-EMR-I.

\bibliographystyle{amsplain}

\end{document}